\documentclass[10pt]{ijnam}
\def\currentvolume{1}
\def\currentissue{1}
\def\currentyear{2009}
\def\month{June}
\pagespan{1}{11}
\copyrightinfo{2009}{} 

\usepackage{mydef}

\begin{document}

\title[Extraction method for Stokes Flow with jumps in the pressure]
{Extraction method for Stokes Flow with jumps in the pressure}

\author[K.S. Chang \and D.Y. Kwak]{Kwang Sung Chang \and Do Young Kwak  }
\address{
  Department of Mathematical Sciences,
  Korean Society for Industrial and Applied Mathematics(KAIST),
  335 Gwahangno (373-1 Guseong-dong), Yuseong-gu, Daejeon 305-701, Republic of Korea
}
\email{ckslove@kaist.ac.kr \and kdy@kaist.ac.kr}
\urladdr{http://mathsci.kaist.ac.kr/$\sim$dykwak/}





\thanks{This research was supported by KOSEF(contract number R01-2007-000-10062-0), Korea.}

\subjclass[2000]{65Z05, 76D07, 76T99}

\abstract{In this paper, we consider a stationary, constant viscosity,
incompressible Stokes flow with singular forces along one or several
interfaces. Assuming only the jumps of the pressure are present along the interface, we develop a new numerical scheme for such a problem. By constructing an approximate singular function and removing it, we can apply a standard finite element method to solve it. A main advantage of our scheme is that one can use a uniform grid. We observe optimal $O(h)$ order for the pressure and $O(h^2)$ order for the velocity.}

\keywords{Stokes equation, singular forces, jumps in the solution, extraction method, discontinuous pressure, uniform grid.}

\maketitle

\section{Introduction}
In recent years, interface problems have become the subject of extensive research \cite{Li1994, Bramble, Chen, Anderson, Jacob, Bunner, TLin, LiYan08, CKW}.
Many interesting physical phenomenons are described by the underlying partial differential equations having interface. For example, when two or more distinct materials or fluid with different conductivities, densities or permeability are involved, model equations often involves discontinuous coefficients to reflect the physical properties \cite{Babuska, unique6, Puckett, Rudman, Renardy, Davalos}. Often, the solutions of these interface problems must satisfy certain interface jump conditions due to physical conservation laws. Many numerical methods to deal with such problems have been proposed \cite{Sethian, Hou_L, Reusken, GWWei}. Because of discontinuity of the solutions, standard numerical method do not yield accurate solutions even when fitted grid are used \cite{unique8, Reusken, GersXfem08}.

In this paper, we propose an accurate, fast finite element method  for Stokes problems with pressure jump conditions across a given interface which divides the domain into two parts. Any regular finite element meshes including uniform meshes are allowed. The idea is to consider certain singular function in a neighborhood of the interface whose jumps match the given jump conditions. By subtracting this function from the variational form, we obtain a new variational problem in which the solution has no jumps. To develop a numerical scheme, we choose one such singular function and construct its approximation. The process consists of two main steps: First, we construct a piecewise linear function satisfying the jump conditions in a small strip near the interface. Then we extend it into whole region in some reasonable way. One of the natural method is to solve a harmonic/biharmonic equation in one of the subdomains. It is quite natural to use finite element methods; However, the usual piecewise linear continuous finite element cannot be used because of discontinuity of the data near the interface. Instead, Crouzeix-Raviart(CR) nonconforming element \cite{Crouzeix, Ciarlet, Girault} which uses the midpoint of each edge as degree of freedom is appropriate in this case (see Section 4.1 for details).

The next step is to subtract it from the variational formulation of Stokes problem which leads to standard variational form of Stokes problem. Some advantages of our scheme are:
\begin{itemize}
\item We can use uniform mesh which is very efficient for moving interface such as time dependent problem.
\item Neither do we need adaptive mesh nor do we need extra degrees of freedom such as XFEM \cite{Reusken}, yet our method achieve $O(h^2)$ for velocity and $O(h)$ for pressure, which is optimal with lowest order finite element.
\item The cost of constructing the discrete singular part is cheap since the work is equivalent to solving a Laplace problem in a subdomain.
\item After subtracting the singular part, the remaining task is equivalent to solving standard Stokes equation, hence it can be incorporated into the existing software.
\end{itemize}

The outline of the paper is as follows: We start the formulation of the problem in Section 2, continued by a brief discussion of the removing singularity in the weak form of the Stokes equation. Section 3,4 describe the extraction method and its numerical scheme. In Section 5, some numerical results are shown. Conclusions follow in Section 6.

\section{Model Stokes problem}
\begin{figure}[t]  \label{domain}
\begin{center}
      \psset{unit=2cm}
      \begin{pspicture}(-1,-1)(1,1)
        \pspolygon(1,1)(-1,1)(-1,-1)(1,-1)
        \pscurve(0.5,0)(0.3,0.3)(-0.6,-0.1)(-0.2,-0.5)(0.52,-0.08)(0.5,0)
        \rput(0,0){\scriptsize$\Omega^-$}
        \rput(-0.3,0.5){\scriptsize$\Omega^+$}
        \rput(0.75,0){\scriptsize$\Gamma$}
      \end{pspicture}
\caption{Sketch of the domain $\Omega$ for the interface problem}
\end{center}
\end{figure}
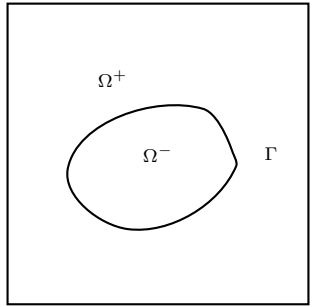
Let $\Omega \subset \RR^2$ be a convex polygonal domain. The domain $\Omega$ is separated into two subdomains
$\Omega^-$ and $\Omega^+$ with $\bar \Omega = \bar \Omega^-
 \cup \bar\Omega^+$ and $\Omega^- \cap \Omega^+ = \emptyset$. We
assume that $\Omega^-$ and $\Omega^+$ are connected and
$\partial\Omega^-\cap \partial\Omega^+ = \emptyset$. The
interface is denoted by $\Gamma = \bar \Omega^- \cap \bar
\Omega^+$. Let %
\begin{equation}
\begin{array}{rcl}
H_0^1(\Omega)&:=&\Brace{v\in H^1(\Omega):\ v=0 \mbox{ on }\partial\Omega},\\[2mm] %
L_0^2(\Omega)&:=&\Brace{q\in L^2(\Omega):\ \int_\Omega q = 0}.
\end{array}
\end{equation}
Then, we want to find the solution $(\uu, p)\in
H_0^1(\Omega)^2 \times L_0^2(\Omega)$ of the stationary homogeneous
Stokes problem for an incompressible viscous fluid confined in
$\Omega$ satisfies: %
\begin{subequations}\label{eq:BasicStokes}
\begin{alignat}{4}
-\mu\Delta \uu + \nabla p &\ =\ &\Bsym g + \Bsym F, &\qquad \mbox{in } \Omega,\label{eq:BasicStokes1}\\ %
\nabla \cdot \uu &\ =\ &0,&\qquad\mbox{in }\Omega,\\ %
\uu&\ =\ &\Bsym 0,&\qquad\mbox{on }\partial\Omega.
\end{alignat}
\end{subequations}%
with $p=p(\xx)$ the pressure, $\uu=\uu(\xx)$ the velocity, $\mu$
the constant viscosity, $\Bsym g\in
L^2(\Omega)$ and $\Bsym F\in H^{-1}(\Omega)$, the external singular force. By \cite{unique1},\hspace{-1mm} \cite{unique2},
this problem has a unique solution $(\Bsym u, p)\in
H_0^1(\Omega)^2\times
L_0^2(\Omega)$.
The external singular force can be written as
$$\Bsym F\equiv (F_1(x,y),F_2(x,y))=\int_\Gamma \Bsym f(s) \delta (\xx -
\XX(s)) \ds,$$ where $\XX(s)$ denotes the interface parameterized by $s$, $\Bsym f(s)$ is the force strength at this point, and $\delta$ is the two-dimensional delta function. In general, this singular force leads to the jumps of the pressure and the velocity. However, dealing with those jumps for both the velocity and pressure is a heavy task, no one seems to have resolved it completely yet. Hence, in this paper, we restrict out attention to a rather simple case where the jumps are restricted the pressure only. So, we assume that the jumps of the pressure on the interface $\Gamma$ are given by
\begin{equation}
\Bracket{p}_\Gamma = J_1(\xx),\quad \Bracket{\pdiff{p}{\nn}}= J_2(\xx).
\end{equation}
Now, we define the related (affine) spaces. Let $D^+:=D\cap\Omega^+,\ D^- :=\ \ D\cap\Omega^-$ and define
$$
\begin{array}{rl}
H_P^m(D):=&H^m(D^+)\cap H^m(D^-),\\[2mm] %
\mathcal{P}^{m;(\gamma_1,\gamma_2)}(D):=&\dis\Brace{p\in H_P^m(D)\Big|\ [p]=\gamma_1,\ \left[\frac{\partial p}{\partial \bold n} \right]=\gamma_2 \mbox{ on } D\cap\Gamma }, \\[2mm] %
\mathcal{P}_0^{m;(\gamma_1,\gamma_2)}(D) :=&\Brace{p\in \mathcal{P}^{m;(\gamma_1,\gamma_2)}(D)|\ p=0 \mbox{ on } \partial D }, %
\end{array}
$$ %
with a piecewise norm
$$
\begin{array}{rcl}
\|p\|^2_{H_P^m(D)}&:=&\|p\|^2_{H^m(D^+)}+\|p\|^2_{H^m(D^-)}.
\end{array}
$$
Here, for any domain in $D$, $H^m (D)$ is the usual Sobolev space of order $m$. Considering the jump conditions, we decompose $p$ as %
\begin{equation}
p=p^0 + p^*
\end{equation}
where $p^0\in H^1(\Omega)\cap L_0^2 (\Omega)$ and $p^*\in\mathcal{P}^{1;(J_1,J_2)}(\Omega)$. In other words, $p$ is splitted into regular part $p_0$ and singular part $p^*$.
Using Green's theorem on each subdomains $\Omega^\pm$, we get
\begin{equation}
\int_{\Omega^\pm} \nabla p \cdot\vv \dxb = - \int_{\Omega^\pm} p \Div\vv \dxb + \int_{\partial\Omega^\pm}p\vv\cdot\nn \ds.
\end{equation}
Multiply a test function $\vv$ in (\ref{eq:BasicStokes1}) and integrating by part in each subdomain $\Omega^+$ and $\Omega^-$, we obtain a weak formulation of our problem as follows: find $(\uu, p^0)\in
H^1_0 (\Omega)^2 \times L_0^2 (\Omega)$ such that %
\begin{subequations}\label{eq:WeakStokes}
\begin{alignat}{4}
a(\uu,\vv) + b(\vv,p^0) &\ =\ (\Bsym g,\vv) - b(\vv,p^*) - <J_1,\vv\cdot\nn>_\Gamma, &\qquad ^\forall \vv \in H_0^1(\Omega)^2,\label{eq:WeakStokes1}\\ %
b(\uu,q) &\ =\ 0, &\qquad ^\forall q \in L^2(\Omega),\label{eq:WeakStokes2}
\end{alignat}
\end{subequations}
where
$$\begin{array}{c}
\dis a(\uu, \vv):=\int_\Omega \mu \nabla\uu : \nabla\vv \dxb, \quad b(\vv, q) := - \int_\Omega q \Div \vv \dxb,\\ %
\dis (\nabla\uu : \nabla\vv) := {\rm tr (\nabla\uu \nabla\vv)}.%
\end{array}$$
The resulting equation is the variational form of a standard Stokes equation with modified right hand side. The problem now is how to find an approximation to $p^*$ and approximate variational form. We will answer this question in the next section.

\section{Extraction Method}
Taking the divergence of the equation (\ref{eq:BasicStokes1}) in each subdomain, we obtain the Poisson equation for the pressure:
\begin{equation}
\nabla^2 (p^0 + p^*)=\Delta\parenthesis{p^0+p^*}=\Div\Bsym g.
\end{equation}
We may split this equation into two equations
\begin{subequations}\label{eq:extraction}
\begin{alignat}{8}
&\qquad &\Delta p^0&= G_2 &\quad &\mbox{ in } \Omega,                    &\qquad     &p^0 \in H^1(\Omega)\cap L_0^2(\Omega),\label{eq:extraction1}\\
&       &\Delta p^*&= G_1 &      &\mbox{ in } \Omega^-\cup\Omega^+,      &           &p^*\in \mathcal{P}_0^{1;(J_1,J_2)}(\Omega),\label{eq:extraction2}
\end{alignat}
\end{subequations}
for some $G_1$ and $G_2$. Since $p^*$ has jump across $\Gamma$, the equations (\ref{eq:extraction2}) hold in $\Omega^+$ and $\Omega^-$ respectively, while (\ref{eq:extraction1}) holds in $\Omega$. But solving this system numerically is not an easy task, since such a splitting is not unique and furthermore, $p^*$ has two interface conditions. In this paper, we propose the following method: First, we consider a narrow strip contained in $\Omega^-$ whose outer boundary coincides with $\Gamma$ (see Figure \ref{fig:subdomains_conti}). Choose any function satisfying jump conditions and restrict it to $S$(call it $p^E$. Then we extend it into the whole $\Omega^-$ by solving the following equation
\begin{equation}
\begin{array}{rcll}
\Delta p^*&=& 0 &\mbox{ in } \Omega^-\setminus S=:\Omega^*,\\
p^*&=& p^E &\mbox{ on } \partial\Omega^*.
\end{array}
\end{equation}
Finally, set $p^*=0$ on $\Omega^+$. We call this scheme an \textit{Extraction Method(EM)}. The remaining task is to find a finite dimensional approximation to $p^*$.

\section{Numerical Scheme}
\begin{figure}
    \subfigure[Subdomains $\Omega^*, S$ in $\Omega^-$.\label{fig:subdomains_conti}]{
        \centering
        \psset {linewidth=0.1pt, yunit=2, xunit=2}
        \begin{pspicture}(-1,-1)(1,1)
        \psframe[fillstyle=solid, fillcolor=white](-1,-1)(1,1) \rput(-0.8,0.8){\footnotesize$\Omega^+$} \rput(0.65,0){\footnotesize$\Gamma$}%
        \pscircle[fillstyle=solid, fillcolor=lightgray, linecolor=blue, linewidth=0.5pt, linestyle=solid](0,0){1.2} \rput(-0.45,0){\footnotesize$S$}  %
        \pscircle[fillstyle=solid, fillcolor=gray, linecolor=gray, linewidth=0.5pt, linestyle=solid](0,0){0.7} \rput(0.05,0){\footnotesize$\Omega^*$} %
        \end{pspicture}
    }
\qquad\qquad
    \subfigure[$\Omega_h^I, \Omega_h^*, \Gamma$ in triangular grids.\label{fig:subdomains_discrete}]{
        \centering
        \psset {linewidth=0.1pt, yunit=2, xunit=2}
        \begin{pspicture}(-1,-1)(1,1)
        \pspolygon[fillstyle=solid, fillcolor=lightgray, linewidth=0.1pt, linecolor=lightgray]%
        (-1,0.5)(-0.5,0.5)(-0.5,1)(0,1)(1,0)(1,-0.5)(0.5,-0.5)(0.5,-1)(0,-1)(-1,0) %
        \pspolygon[fillstyle=solid, fillcolor=gray, linewidth=0.1pt, linecolor=blue](0,0)(0,0.5)(0.5,0) %
        \pspolygon[fillstyle=solid, fillcolor=gray, linewidth=0.1pt, linecolor=blue](0,0)(0,-0.5)(-0.5,0) %
        \multido{\R=-1+0.5}{5}{\psline(\R,-1)(\R,1)}%
        \multido{\R=-1+0.5}{5}{\psline(-1,\R)(1,\R)}%
        \multido{\R=-1+0.5}{5}{\psline(-1,\R)(\R,-1)} %
        \multido{\R=-1+0.5}{4}{\psline(\R,1)(1,\R)} %
        \pscircle[linecolor=blue, linewidth=0.5pt, linestyle=solid](0,0){1.2} %
        \psline[linestyle=dashed]{->}(-1.1,0.1)(-0.8,0.1) \rput(-1.2,0.1){\tiny$\Omega_h^I$} %
        \psline[linestyle=dashed]{->}(1.1,0.8)(0.1,0.2) \rput(1.2,0.8){\tiny$\Omega_h^*$} %
        \psline[linestyle=dashed]{->}(1.1,-0.8)(0.55,-0.25) \rput(1.2,-0.8){\tiny$\Gamma$} %
        \end{pspicture}
    }
    \caption{Some subdomains in Extraction method}\label{fig:subdomains}
\end{figure}
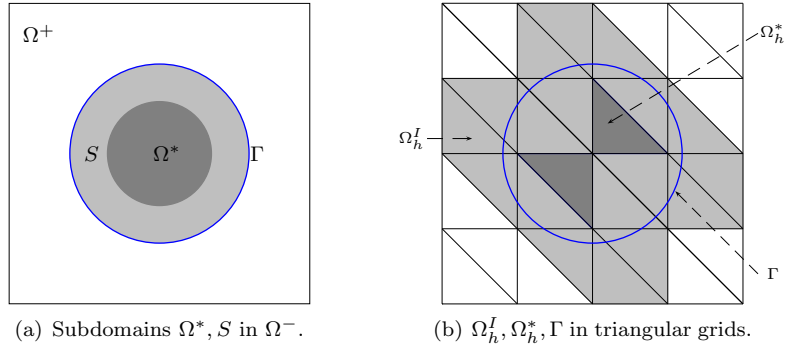
\begin{figure}[t]
\centering
    \subfigure[The numerical interface $\Gamma_h$\label{fig:interface_h}]
    {
        \centering
        \includegraphics[width=6cm]{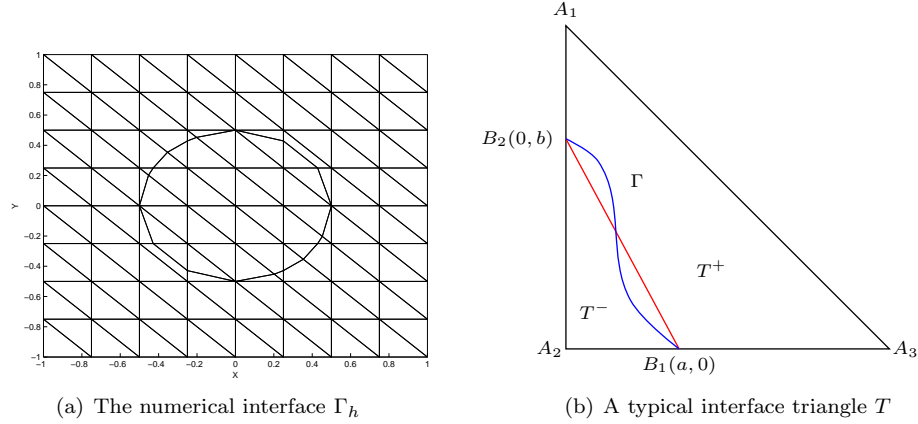}
    }
\qquad\qquad
    \subfigure[A typical interface triangle $T$\label{fig:inter_tri}]
    {
        \centering
        \psset{unit=4.3cm}
        \begin{pspicture}(0,-0.1)(1,1)
          \psset{linewidth=0.5pt}
          \pspolygon(0,0)(1,0)(0,1)
          \psline[linecolor=red](0,0.65)(0.35,0)
          \pscurve[linecolor=blue](0,0.65)(0.1,0.58)(0.2,0.15)(0.35,0)
          \rput(0,1.05){\scriptsize$A_1$}
          \rput(-0.05,0){\scriptsize$A_2$}
          \rput(1.05,0){\scriptsize$A_3$}
          \rput(-0.15,0.65){\scriptsize$B_2(0,b)$}
          \rput(0.09,0.12){\scriptsize$T^-$}
          \rput(0.45,0.25){\scriptsize$T^+$}
          \rput(0.35,-0.05){\scriptsize$B_1(a,0)$}
          \rput(0.22,0.52){\scriptsize$\Gamma$}
        \end{pspicture}
    }
    \caption{Interface elements in the triangular mesh.}\label{fig:interface_elements}
\end{figure}
Let $\mathcal{T}_h$ be the usual quasi-uniform finite element triangulations of the domain $\Omega$. For any element $T\in\mathcal{T}_h$, we call an element
$T$ an \textit{interface element} if the interface $\Gamma$ passes through the interior of $T$, otherwise we call $T$ a \textit{non-interface element} and we call an edge $e\in\partial T$ an \textit{interface edge} if the interface $\Gamma$ passes through the interior of $e$, otherwise we call $e$ a \textit{non-interface edge}. Now, we introduce some notations:
$$\begin{array}{rcl}
\partial^I T&=&\mbox{the set of all interface edges of an element } T, \\[2mm] %
\mathcal{T}^I_h&=&\mbox{the set of all interface elements}, \\[2mm] %
\mathcal{T}^N_h&=&\mbox{the set of all non-interface elements}, \\[2mm] %
\Omega_h^I&=&\dis\bigcup_{T\in\mathcal{T}_h^I} T,\\  %
\Omega_h^-&=&\dis\bigcup_{T\in\mathcal{T}_h} T^-,\\   %
\Omega_h^*&=&\dis\Omega_h^- \setminus \Omega_h^I.   %
\end{array}$$
Even though the interface $\Gamma$ is a curve in general, we replace the part of interface in $T$ by the line segment connecting the intersection points with $\partial T$. Therefore, the interface $\Gamma$ is replaced by
its polygonal approximation $\Gamma_h$.
Henceforth, $\Gamma$ is always assumed to be $\Gamma_h$(see Figure \ref{fig:interface_h}).

\subsection{Construction of $p_h^*$}

We construct $p_h^*$ in two steps. First, we will consider $p_h^*$ in $\Omega_h^I$.
Suppose $T$ is an interface element. For simplicity, we assume the
three vertices are given by $A_1=(0,1)$, $A_2=(0,0)$, $A_3=(1,0)$
(see Figure \ref{fig:inter_tri}). For any element $T$ in general
position, all the constructions to be presented below carries over
 through affine equivalence. Assume the interface meets with the
element's edges at points $B_1$ and $B_2$.

Let $\ell_i$ be the usual linear Lagrange nodal basis function
associated with the vertex $A_i$ for $i=1,2,3$. Then $\ell_1=y $,
$\ell_2=1-x-y$, $\ell_3=x$. In the domain $\Omega_h^I$, we write $\psi^*$ as the following form: %
\begin{equation}\label{def:basis}
\psi^* = \left\{
\begin{array}{rcll}
    \psi^{*-} &=& \alpha_1\ell_1+ \alpha_2\ell_2 + \alpha_3\ell_3  &\text{in $T^-$,}\\
    \psi^{*+} &=& 0 &\text{in $T^+$}. %
\end{array}\right.
\end{equation}
Now imposing the jump conditions on $\Gamma$, we have:
\begin{subequations}\label{condi_impose}
\begin{alignat}{2}
\psi^{*-}(B_1)-\psi^{*+}(B_1) &= J_1(B_1),\label{condi_impose1}\\
\psi^{*-}(B_2)-\psi^{*+}(B_2) &= J_1(B_2),\label{condi_impose2}\\ %
\quad \pd{\psi^{*-}}{\nn}-\pd{\psi^{*+}}{\nn} \ \ &=\ J_2(B_0),\label{condi_impose3} %
\end{alignat}
\end{subequations}
where $B_0$ is the midpoint of $\overline{B_1 B_2}$. Then we have three unknowns $\alpha_i, \ i=1,2,3$ in three equations  (\ref{condi_impose}a,b,c).
Thus we can find coefficients $\Brace{{\alpha}_{i=1,2,3}}$. Note that two end points of the line segment $\Gamma_{seg}$ are
located on the interface $\Gamma$, and hence the interface condition
$[p](\xx)=J_1(\xx)$ is enforced exactly at these two end points,
i.e., the point jump conditions (\ref{condi_impose}a,b) gives
\begin{eqnarray*}
(1-a)\alpha_2 + a\alpha_3 &=& J_1(B_1),\\
b\alpha_1 + (1-b)\alpha_2 &=& J_1(B_2).
\end{eqnarray*}
The second condition on the interface segment $\Gamma_{seg}$ is the
flux continuity. Hence, the derivative jump condition (\ref{condi_impose3})
becomes
$$
(\alpha_1\Grad\ell_1+\alpha_2\Grad\ell_2+\alpha_3\Grad\ell_3)\cdot\mathbf{n}_{\overline{\textrm{\tiny{$B_1B_2$}}}}
=J_2(B_0).
$$
Since these conditions are represented by the following matrix equation:
\begin{eqnarray} \label{eq:local_basis}
\left[\begin{array}{ccc}
0 & 1-a & a \\
b & 1-b & 0 \\
a & -a-b & b \\
\end{array}\right]
\left(\begin{array}{c}
\alpha_1 \\ \alpha_2 \\ \alpha_3 \\
\end{array}\right)
= \left(\begin{array}{c}
J_1(B_1) \\ J_1(B_2) \\ J_2(B_0) \\
\end{array}\right)
=: \left(\begin{array}{c}
\jmath_1 \\ \jmath_2 \\ \wp \\
\end{array}\right),
\end{eqnarray}
 the coefficient $\alpha_i$'s of $\psi^*$ are determined
by the following formula:
\begin{equation}
\psi^*: \left\{\begin{array}{rcl}
\alpha_1&=&\dis \parenthesis{b(b-1)\jmath_1 + (a^2+b)\jmath_2-a(-b-1)\wp}/K \\[2mm] %
\alpha_2&=&\dis \parenthesis{b^2\jmath_1  + a (a \jmath_2-b \wp)}/K \\[2mm] %
\alpha_3&=&\dis \parenthesis{(a+b^2)\jmath_1+(a-1)(a\jmath_2-b\wp)}/K  %
\end{array}\right.
\end{equation}
where $K=a^2+b^2$. We do this for every $T\in\Omega_h^I$. Having constructed $p_h^E$ on $\Omega_h^I$ (lightly shaded region in Figure \ref{fig:subdomains_discrete}), we now need to extend it into $\Omega_h^*$ (dark shaded region in Figure \ref{fig:subdomains_discrete}), which will be done by solving the Laplace equation
\begin{equation}\label{eq:Lapalce_pressure}
\begin{array}{rcll}
\Delta p^*&=& 0 &\mbox{ in } \Omega_h^*,\\
p^*&=& p_h^E &\mbox{ on } \partial\Omega_h^*,
\end{array}
\end{equation}
numerically. Numerical methods for this problem are well known; However, due to the discontinuity in the boundary data $p_h^E$ on $\partial\Omega_h^*$, the usual nodal finite element space cannot be used. Instead, one can use Crouzeix-Raviart $P_1$-nonconforming finite element space\cite{unique1} where the linear basis function has the degrees of freedom at the midpoints of edges. We denote it by $\dot P_h$. The finite element solution of (\ref{eq:Lapalce_pressure}) on $\Omega_h^*$ is denoted by $p_h^*$. Together with the construction above on $\Omega_h^I$, we have obtained $p_h^*$ in $\Omega_h^-$.

\subsection{Variational form after removing $p_h^*$}

In this section, we explain how to remove the discrete singular part $p_h^*$ from the weak formulation (\ref{eq:WeakStokes}). First of all, we define some discrete spaces:
$$
\begin{array}{rcl}
P_h^0&:=& \dis \Brace{p\in L_2(\Omega_h)\vert\ p|_T \mbox{ is constant for every } T\in T_h},\\
\what P_h^I&:=& \dis \Brace{p\in L_2(\Omega_h^I)\vert\ p|_{T^\pm} \mbox{ is linear for every } T\in T_h^I},\\[2mm]
\dot P_h(D)&:=& \left\{p\in L_2(D) \vert\ p \mbox{ is linear for every } T\cap D,\ T\in T_h\right.,\\
&& \hspace{2cm}\left. p \mbox{ is continuous at the midpoints of the triangle edges}\right\},\\
P_h^E &:=& \what P_h^I \oplus \dot P_h(\Omega_h^*),\\
\Bsym V_h &:=& (\dot P_h(\Omega_h))^2.
\end{array}
$$
We now replace $p^*$ in (\ref{eq:WeakStokes}) by $p_h^*$, but not without caution: Since $p_h^*$ is now discontinuous along edges of $T$, we include line integrals. Thus we replace (\ref{eq:WeakStokes1}) by the following form
\begin{equation}
a(\uu,\vv) + b(\vv,p^0) \ =\ (\Bsym g,\vv) - b(\vv,p_h^*) + \mathcal{J}(\vv,p_h^*),
\end{equation}
where
$$
\mathcal{J}(\vv,p_h^*):= - <\Bracket{p_h^*},\vv\cdot\nn>_\Gamma + \sum_{T\in\TT_h^I}\int_{\partial^I T} p_h^* \vv\cdot\nn \ds.
$$
The resulting equation can be solved again by a standard finite element method for Stokes problem: The simplest and most natural finite element method in this setting is the Crouzeix-Raviart finite element pair $\Bsym V_h \times P_h^0$. Thus, we have the following discrete Stokes problem:
Find $(\uu_h, p_h^0)\in
\Bsym V_h \times P_h^0$ such that %
\begin{subequations}
\begin{alignat}{4}
a_h(\uu_h,\vv_h) + b_h(\vv_h,p_h^0) &\ =\ (\Bsym g,\vv_h) - b_h(\vv_h,p_h^*) + \mathcal{J}(\vv_h,p_h^*), &\quad^\forall \vv_h \in \Bsym V_h,\\ %
b_h(\uu_h,q_h) &\ =\ 0, &\quad^\forall q_h \in P_h^0,
\end{alignat}
\end{subequations}
where
$$
\begin{array}{c}
\dis a_h(\uu_h, \vv_h):=\sum_{T\in\TT}\int_T \mu \nabla\uu_h : \nabla\vv_h \dxb, \quad b_h(\vv_h, q_h) := - \sum_{T\in\TT}\int_{T} q_h \Div \vv_h \dxb. %
\end{array}
$$

\section{Numerical experiments}
In all of the experiments, the domain is a square and triangularized by uniform triangle grids with $h_x=h_y=1/2^{n-1}$ for $n=3,\cdots,8$. In order to describe the interface, we consider the level-set function $\Phi(\xx)$ for the interface
$\Gamma$ which is assumed to be smooth. Let $\Phi: \Omega
\rightarrow \RR$ be a continuous function such that
\begin{equation}
\Phi(\xx)=\left\{ %
\begin{array}{ll}
<0&\quad \xx \mbox{ in } \Omega^-,\\[2mm] %
=0&\quad \xx \mbox{ on } \Gamma, \\[2mm] %
>0&\quad \xx \mbox{ in } \Omega^+. %
\end{array} %
\right.
\end{equation}
We assume that $\Phi(\xx)$ is smooth and $\nabla\Phi$ is not zero in
any neighborhood of the interface $\Gamma$. Then the unit normal
vector $\nn(\xx)$ is represented by
$\frac{\nabla\Phi}{|\nabla\Phi|}$. The experiments in this subsection show that the method is robust.
\begin{example}[Constant Jump] The level-set function $\Phi(\xx)$, the jumps of the pressure and the boundary condition of the velocity are given as follows:
\begin{equation*}
\begin{array}{rll}
    \Phi(\xx) &=&\sqrt{(x-0.5)^2+(y-0.5)^2}-0.25,\\[2mm]
    p&=&
    \left\{
    \begin{array}{ll}
        \dis 150(x-1/2)(y-1/2) ,& \mbox{in } \Omega^+\\[2mm]
        \dis 150(x-1/2)(y-1/2) + 30,& \mbox{in } \Omega^-
    \end{array}\right.\\[5mm]
    \uu&=&
    \left\{
    \begin{array}{l}
        u_1 = -256x^2(x-1)^2y(y-1)(2y-1),\\[2mm]
        u_2 = 256y^2(y-1)^2x(x-1)(2x-1),
    \end{array}
    \right.\mbox{in } \Omega^\pm\\[5mm]
    \uu&=&{\Bsym 0} \quad \mbox{ on } \partial\Omega.
\end{array}
\end{equation*}
with the domain
$\Omega=\Bracket{0,1}\times\Bracket{0,1}$.
\end{example}
We observe the robust first order for the pressure and second order convergence for the velocity with $L^2$-norm.

\begin{table}[h]
\footnotesize
\caption{(Constant jump) 1st order for the pressure and 2nd order for the velocity with $\norm{\cdot}_{L^2}$} %
\centering%
{
\begin{tabular}{c cc cc cc cc cc}\toprule
&$N_x\times N_y$    && $||p-p_h||_{L^2}$          && Order   && $||u-u_h||_{L^2}$         && Order &\\ \midrule %
&$8\times8$         && $5.1852 \times 10^{-0}$    && -       && $1.6315 \times 10^{-1}$   && -     &\\ %
&$16\times16$       && $2.0473 \times 10^{-0}$    && 1.34    && $4.7486 \times 10^{-2}$   && 1.78  &\\ %
&$32\times32$       && $8.1309 \times 10^{-1}$    && 1.33    && $1.2568 \times 10^{-2}$   && 1.92  &\\ %
&$64\times64$       && $3.6409 \times 10^{-1}$    && 1.16    && $3.2055 \times 10^{-3}$   && 1.97  &\\ %
&$128\times128$     && $1.7276 \times 10^{-1}$    && 1.07    && $8.0694 \times 10^{-4}$   && 1.99  &\\ %
&$256\times256$     && $8.6181 \times 10^{-2}$    && 1.00    && $2.0220 \times 10^{-4}$   && 2.00  &\\ %
\bottomrule %
\end{tabular}
}
\end{table}

\begin{figure}[h]
\centering
\subfigure[$p_h^0$ in $\Omega$]{\includegraphics[width=6cm]{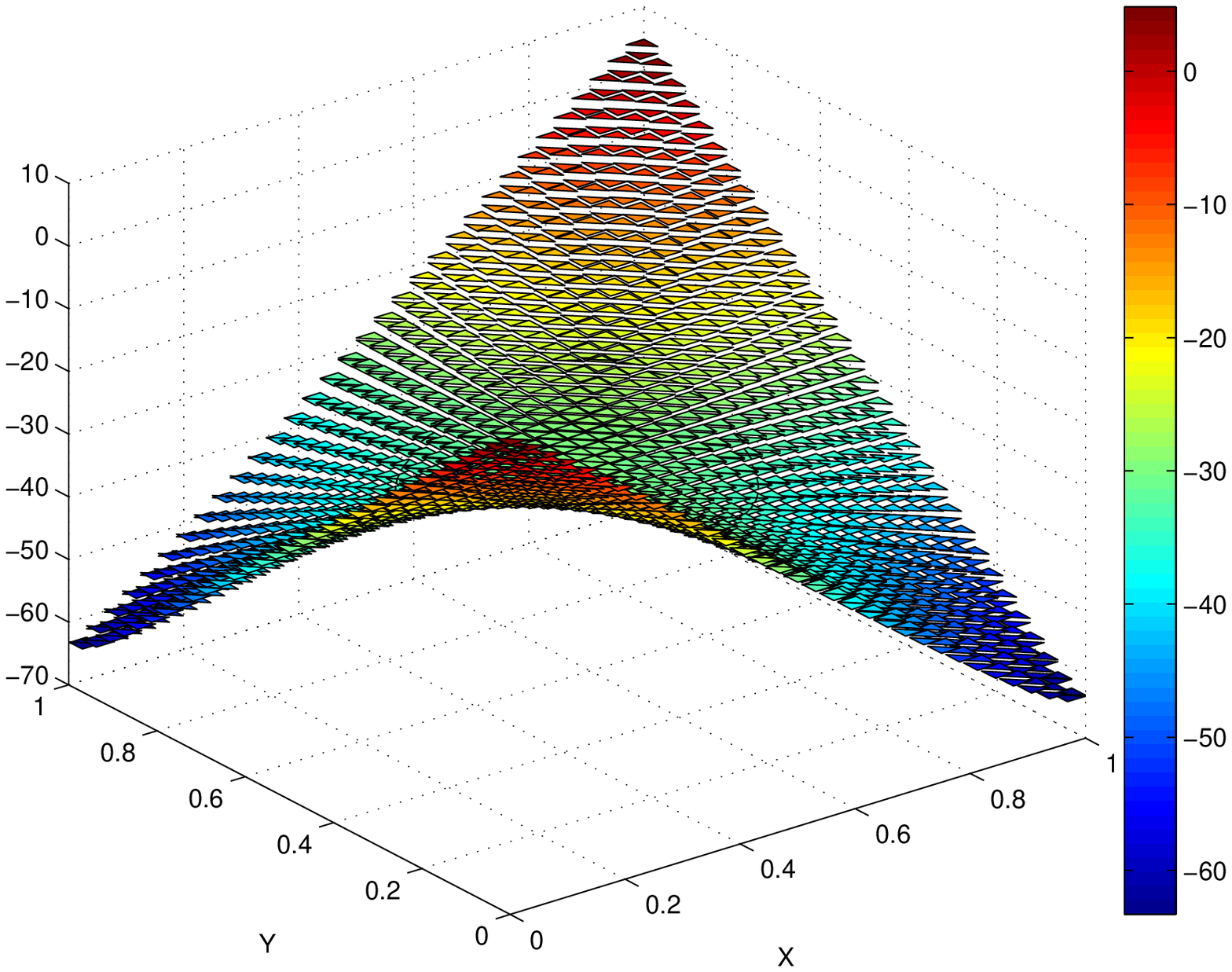}}
\subfigure[$p_h^*$ in $\Omega_h^I$]{\includegraphics[width=6cm]{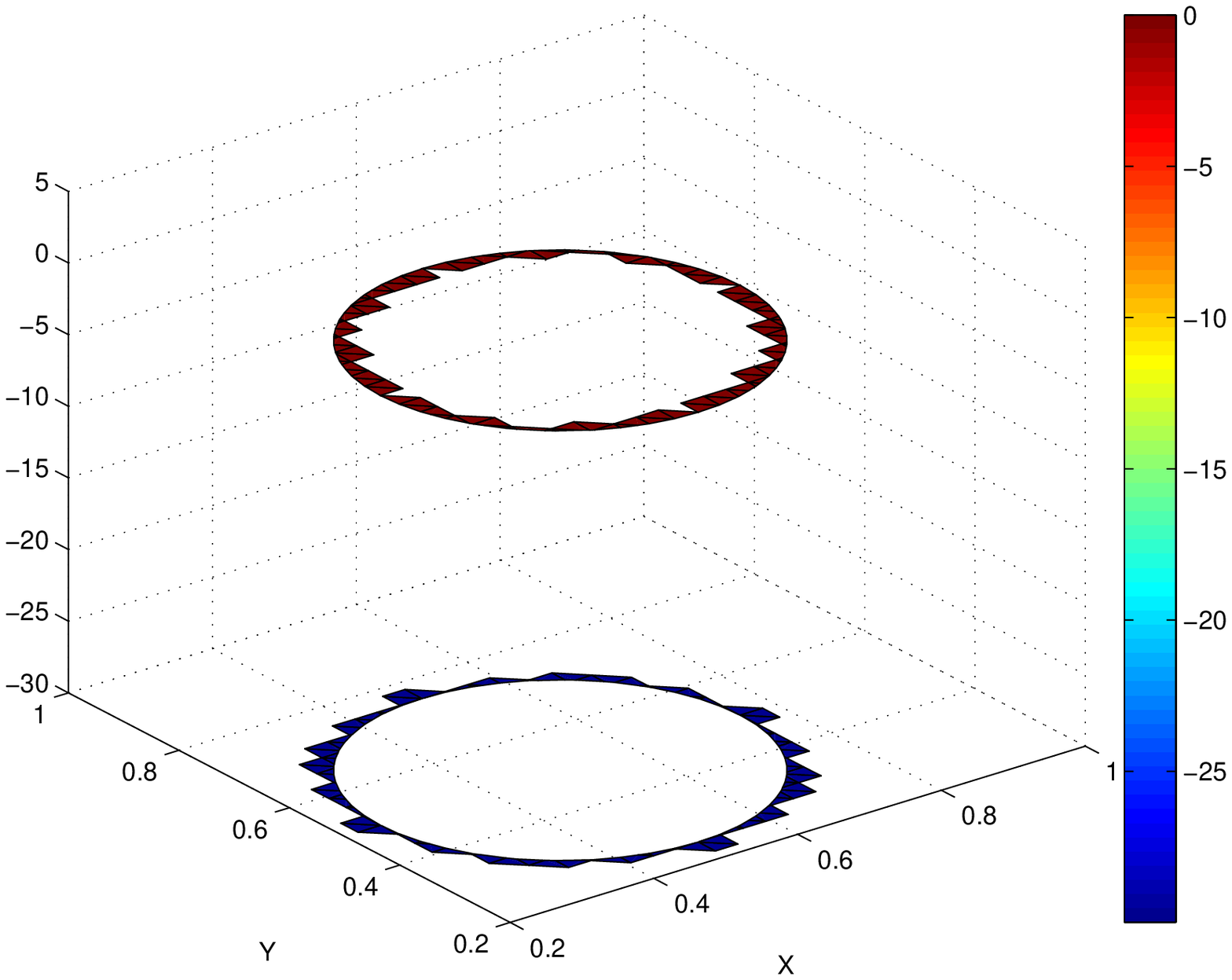}}
\subfigure[$p_h^*$ in $\Omega$]{\includegraphics[width=6cm]{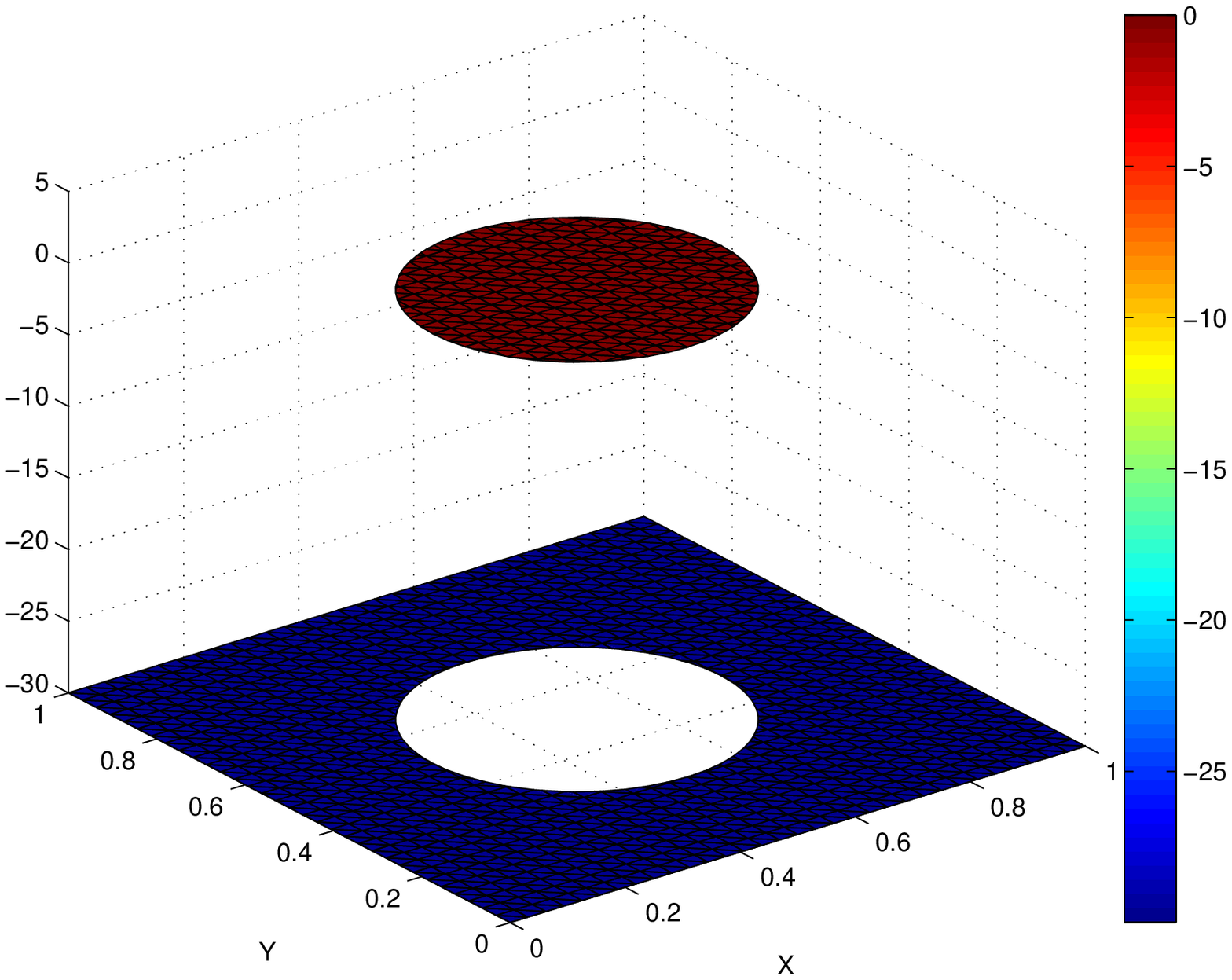}}
\subfigure[$p_h = p_h^0 + p_h^*$ in $\Omega$]{\includegraphics[width=6cm]{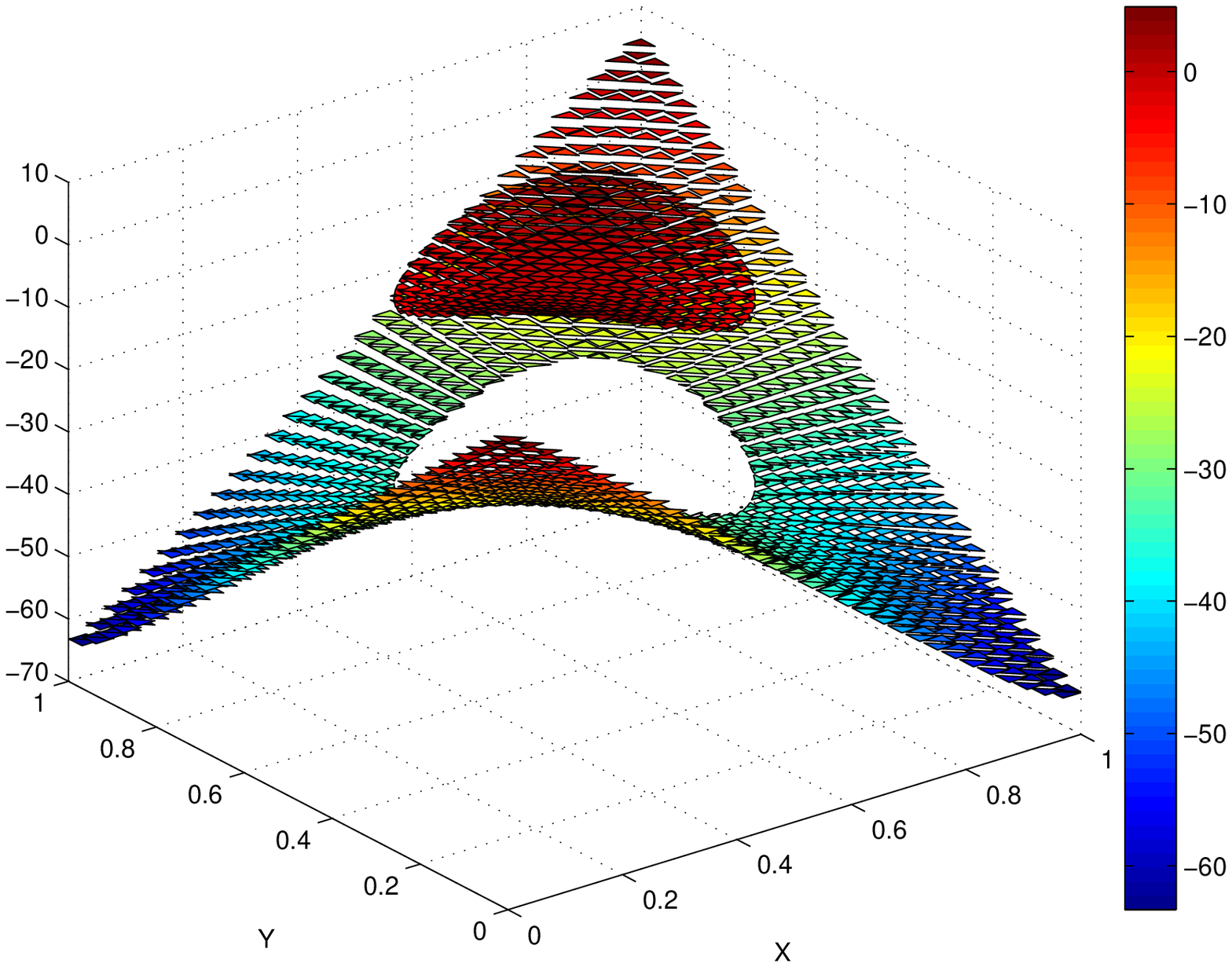}}
\caption{$p_h^0,\ p_h^*,\ p_h$ in the constant jump case.}
\end{figure}

\begin{example}[Noncontant Jump] The level-set function $\Phi(\xx)$, the jumps of the pressure and the boundary condition of the velocity are given as follows:
\begin{equation*}
\begin{array}{rll}
    \Phi(\xx) &=&\sqrt{x^2+y^2}-0.5,\\[2mm]
    \Bracket{p}_\Gamma&=&20\sin(x^2y)-x^2y,\\[2mm]
    \dis\Bracket{\pdiff{p}{\nn}}_\Gamma&=&3x^2y(20\cos(x^2y)-1),\\[2mm]
    \Bsym g&=&
    \left\{
    \begin{array}{ll}
        \dis 20x\cos(x^2y)(2y,x),& \mbox{in } \Omega^-\\[2mm]
        \dis x(2y,x),& \mbox{in } \Omega^+
    \end{array}
    \right.,\\[4mm]
    \uu&=&{\Bsym 0} \quad \mbox{ on } \partial\Omega.
\end{array}
\end{equation*}
with the domain
$\Omega=\Bracket{-1,1}\times\Bracket{-1,1}$.
\end{example}
We again have similar optimal convergence behavior.

\begin{table}[h]
\footnotesize
\caption{(Nononstant jump) 1st order for the pressure and 2nd order for the velocity with $\norm{\cdot}_{L^2}$} %
\centering%
{
\begin{tabular}{c cc cc cc cc cc}\toprule
&$N_x\times N_y$    && $||p-p_h||_{L^2}$          && Order   && $||u-u_h||_{L^2}$         && Order &\\ \midrule %
&$8\times8$         && $2.4718 \times 10^{-1}$    && -       && $1.6473 \times 10^{-2}$   && -     &\\ %
&$16\times16$       && $9.5051 \times 10^{-2}$    && 1.38    && $4.8611 \times 10^{-3}$   && 1.76  &\\ %
&$32\times32$       && $4.2345 \times 10^{-2}$    && 1.17    && $1.5310 \times 10^{-3}$   && 1.67  &\\ %
&$64\times64$       && $2.0367 \times 10^{-2}$    && 1.06    && $4.2185 \times 10^{-4}$   && 1.86  &\\ %
&$128\times128$     && $1.0152 \times 10^{-2}$    && 1.00    && $1.0879 \times 10^{-4}$   && 1.96  &\\ %
&$256\times256$     && $5.0738 \times 10^{-3}$    && 1.00    && $2.7698 \times 10^{-5}$   && 1.97  &\\ %
\bottomrule %
\end{tabular}
}
\end{table}

\begin{figure}[h]
\centering
\subfigure[$p_h^0$ in $\Omega$]{\includegraphics[width=6cm]{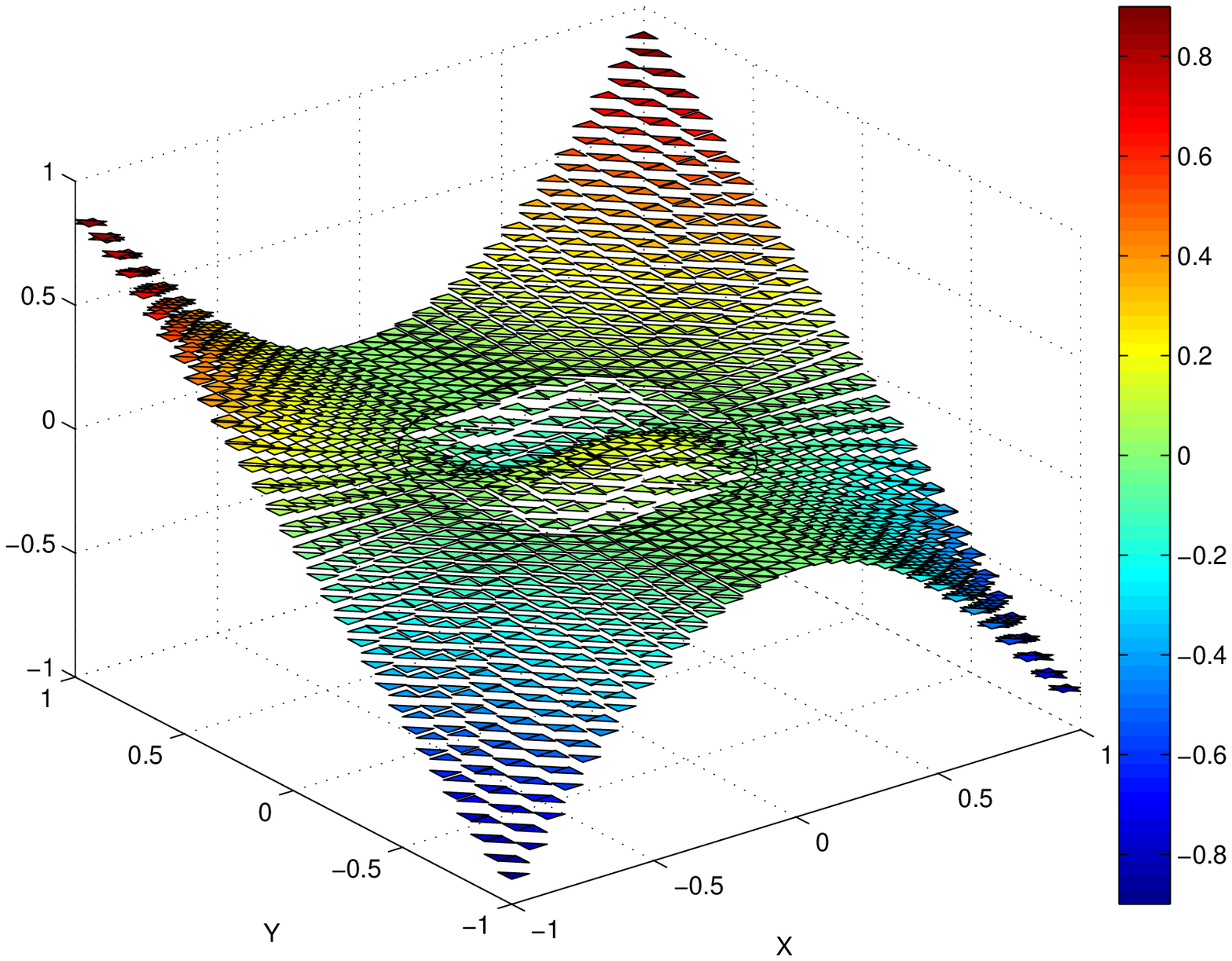}}
\subfigure[$p_h^*$ in $\Omega_h^I$]{\includegraphics[width=6cm]{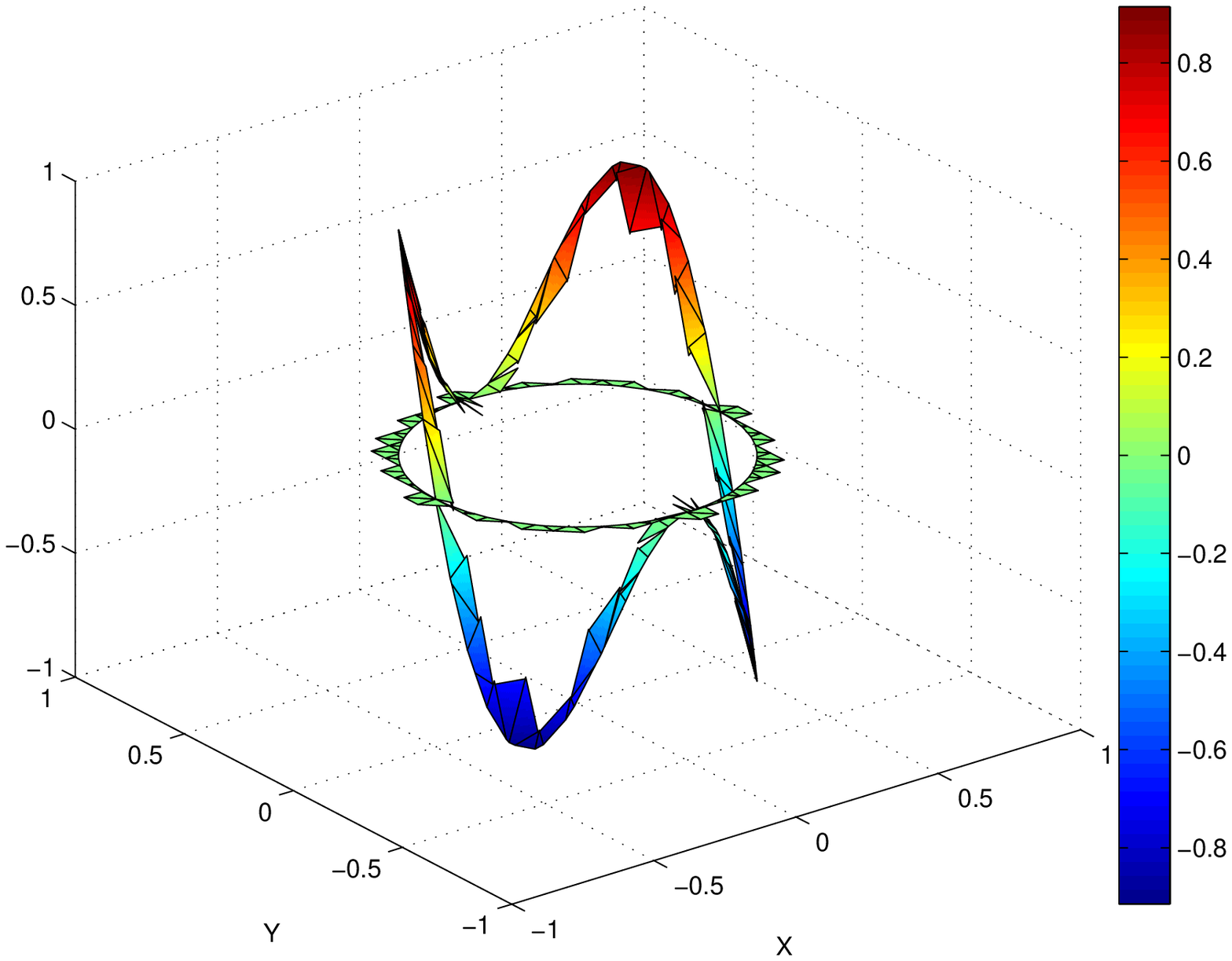}}
\subfigure[$p_h^*$ in $\Omega$]{\includegraphics[width=6cm]{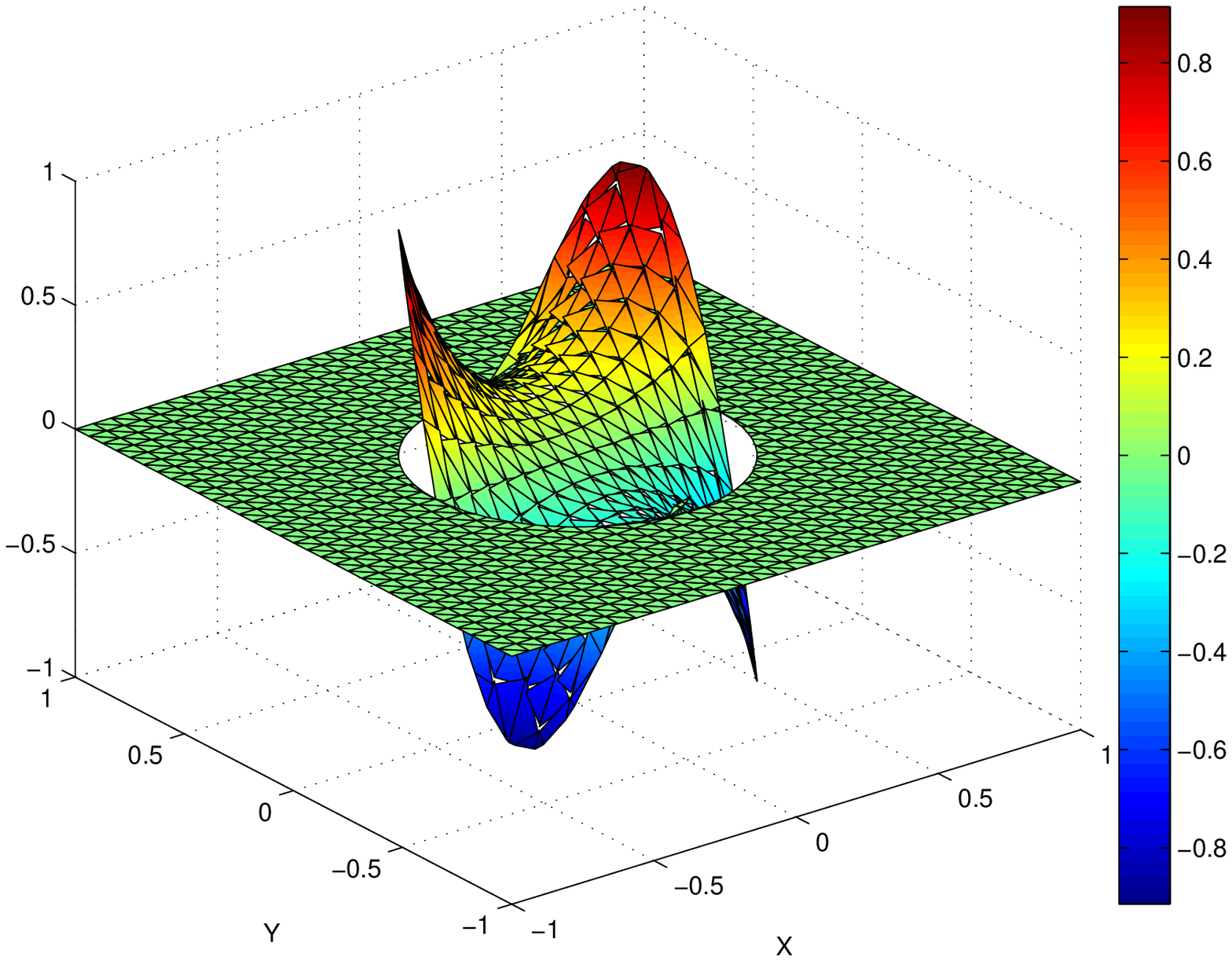}}
\subfigure[$p_h = p_h^0 + p_h^*$ in $\Omega$]{\includegraphics[width=6cm]{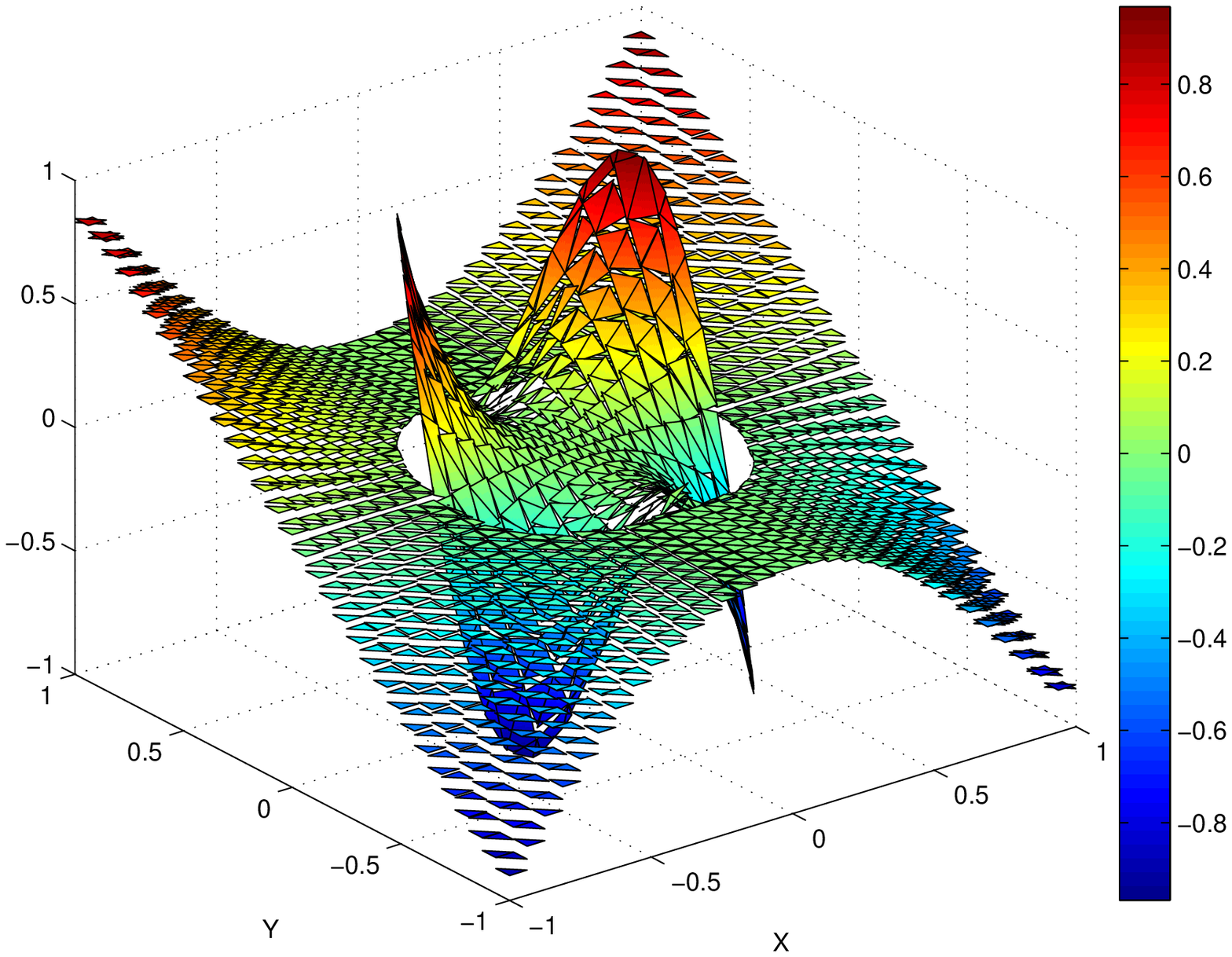}}
\caption{$p_h^0,\ p_h^*,\ p_h$ in the nonconstant jump case.}
\end{figure}

\section{Conclusions} In this paper, we have introduced a new numerical method of solving Stokes interface problems having jumps in the pressure. The first step is to construct a piecewise linear function having small support in $\Omega_h^-$ near the interface which satisfy the jump conditions. The second step is to extend it into $\Omega_h^*$ by solving a discrete Laplace equation with $P_1$-nonconforming finite element. Then removing it from the original variational form, we obtain a Stokes problem with no jumps. The equation is then solved with the Crouzeix-Raviart nonconforming finite element pair. Our scheme is very effective since we can use any shape regular grid, not necessarily fitted grid. We have provided some numerical examples which show the optimal $O(h^2)$ error for velocity and $O(h)$ for pressure.

\end{document}